# On Generalised Inflations of Right Modular Groupoids

## R. A. R. Monzo[1]

**Abstract.** We prove that if a right modular groupoid $U$ is cancellative or has a left identity then a right modular generalised inflation of $U$ is an inflation of $U$.

## 1. Introduction

Since groupoids in the groupoid variety determined by the identity $xy \cdot z = zy \cdot x$ were studied by Kazim and Nazeerudin in 1977 under the name of *left almost semigroups* they have attracted a wide interest [6]. They have appeared under the names of *Abel Grassmann's*, *left invertive* and *right modular groupoids* [1,3-5,7-9,11-17]. The concept of a *generalized inflation* $G$ of its subsemigroup $U$ was introduced in 2000 by Clarke and Monzo as an extension of the concept of an *inflation* of a semigroup [2]. It had been known that a semigroup $S$ has a semilattice of groups as its square if and only if $S$ is an inflation of $S^2$ and $S^2$ is a semilattice of groups [18]. In [2] it was proved that a semigroup $S$ has a square that is a union of groups if and only if $S$ is a (symmetric) generalised inflation of $S^2$ and $S^2$ is a union of groups.

Wang and Wismath pointed out that the concept of a generalised inflation can be extended to any algebra and therefore to groupoids; hence, the resultant generalised inflation need not be a semigroup, even if $U$ is a semigroup [19]. An interesting result that links with those results on semigroups mentioned in the above paragraph was obtained in [11], where it was proved that if a generalised inflation $G$ of its subgroupoid $U$ is a right modular groupoid and $U$ is a union of groups then $G$ is an inflation of $U$.

In this paper we prove that a right modular generalised inflation $G$ of a cancellative groupoid $U$, or of a groupoid $U$ with a left identity element, is an inflation of $U$.

## 2. Notation, definitions and preliminary results

By a groupoid we shall mean a set $G$ with a product $* : G \times G \to G$ and we shall denote $*(x, y)$ by $xy$ or $x \cdot y$. We define $x^1 \equiv x$ and, by induction, $x^n \equiv x \cdot x^{n-1}$ for any $n \in \{2, 3, ...\}$. For example, $x^4 = x \cdot (x \cdot x^2)$. A groupoid $G$ will be called idempotent if $x = x^2$ ($x \in G$). A groupoid $G$ is called a **right modular groupoid** if $G$ satisfies the identity $xy \cdot z = zy \cdot x$.

**Definition.** A groupoid [semigroup] $G$ *is an inflation of its subgroupoid [subsemigroup]* $U$ if $G = \bigcup G_u (u \in U)$, where *(1)* $u \in G_u (u \in U)$, *(2)* $G_u \cap G_v = \varnothing$ if $u \neq v$ and *(3)* $x \in G_u$ and $y \in G_v$ implies $xy = uv$.

**Definition.** A groupoid [semigroup] $G$ *is a generalised inflation of its subgroupoid [subsemigroup]* $U$ if $G = \bigcup G_u (u \in U)$, where *(1)* $u \in G_u (u \in U)$, *(2)* $G_u \cap G_v = \varnothing$ if $u \neq v$, *(3)* for every $x \in G_u$ there exist $\alpha_x$ and $\beta_x$, right and left mappings respectively on $U$, such that for every $x \in G_u$ and $y \in G_v$, $xy = v\alpha_x \cdot \beta_y u$ and **(4)** $u \in U$ implies $\alpha_u = c_u = \beta_u$, where $c_u$ is a constant mapping on $U$ that sends every element to $u$. This implies that a generalised inflation $G$ of its subgroupoid $U$ preserves the product in $U$; that is, for $\{u,v\} \subseteq U$ the product of $u$ and $v$ in $G$ is the product of $u$ and $v$ in U.

[1] Flat 10 Albert Mansions, Crouch Hill, London N8 9RE, United Kingdom ; bobmonzo@talktalk.net

**Definition.** A generalised inflation G of its subgroupoid [subsemigroup] U is called *a symmetric generalised inflation of its subgroupoid [subsemigroup]* U if $u\alpha_x = \beta_x u$ for every $x \in G$ and $u \in U$. If a symmetric generalised inflation G of U satisfies $u\alpha_x = v\alpha_x$ for every $x \in G$ and $\{u, v\} \subseteq U$ then G is *a constant generalised inflation of* U.

Note that a constant generalised inflation G of its subgroupoid U is an inflation of U [2]. Conversely, an inflation G of its subgroupoid U can be expressed as a constant generalised inflation G of U. However there are non-constant generalised inflations G of U that are inflations of U. As stated in the introduction, if a semigroup U is considered as a groupoid then it may have generalised inflations that are ***not*** semigroups.

The following results will be used later in the proof of Theorems 1 and 3:

**Result 1 ([6]).** *If* G *is a right modular groupoid and* $\{a, b, c, d\} \subseteq G$ *then* $ab \cdot cd = ac \cdot bd$.

**Result 2. ([11]).** *If* G *is a constant generalised inflation of its subgroupoid* U *then* G *is an inflation of* U.

**3. Generalised inflations of right modular groupoids**

Although there may be many generalised inflations indeed of a given right modular groupoid, if we require them to also be right modular groupoids then there are far fewer such generalised inflations. In certain cases, the only ones are the inflations, as in the following example of the unique right modular, idempotent groupoid of order 4. It is non-associative and is a quasigroup. Non-associative, right modular, idempotent groupoids of order $\leq 3$ do not exist [17].

**Example 1:** $W = \{1, 2, 3, 4\}$ with multiplication table as follows.

| W | 1 | 2 | 3 | 4 |
|---|---|---|---|---|
| 1 | 1 | 3 | 4 | 2 |
| 2 | 4 | 2 | 1 | 3 |
| 3 | 2 | 4 | 3 | 1 |
| 4 | 3 | 1 | 2 | 4 |

**Table 1:** $W = \{1, 2, 3, 4\}$

If we take $x \notin W$ and let $H = W \cup \{x\}$ then there are over 282,270,000,000,000,000 groupoids on the set H. If we require the multiplication in W to remain intact then
  (1) there are $5^9 = 1,953,125$ different groupoids that can be formed on the set H, some of which may be isomorphic,
  (2) there are $4^9 = 262,144$ different groupoids H in which $H^2 \subseteq W$ and
  (3) there are less than $4^9 + 1$ generalised inflations H of W.

Clearly, an inflation of a right modular groupoid is right modular. Which of the generalised inflations H of W are right modular generalised inflations? Since W is cancellative there are only 4, the inflations of W, as we now prove.

**Theorem 1:** *If* G *is a right modular generalised inflation of its right cancellative subgroupoid* U *then* G *is an inflation of* U.

*Proof:* Let $\{u,v,w\} \subseteq U$ and let $x \in G - U$. Then $ux \cdot v = vx \cdot u$ and so $(u\beta_x u)v = (v\beta_x v)u = (u\beta_x v)v$. Therefore, $u\beta_x u = u\beta_x v$. Also, $xu \cdot v = vu \cdot x$ and so $(u\alpha_x u)v = (vu)\beta_x(vu) = (vu)\beta_x w = [(\beta_x w)u]v$. Therefore, $u\alpha_x = \beta_x w$, for every $\{u,w\} \subseteq U$. Thus, G is a constant generalised inflation of U and, by Result 2, an inflation of U. This completes the proof of Theorem 1.

It is straightforward to show that a right modular groupoid is cancellative if and only if it is right cancellative. Also, the following result has been proved in [11]:

**Theorem 2:** *Let* U *be a right modular groupoid and a union of groups. Then any right modular generalised inflation of* U *is an inflation of* U. *In fact, any right modular groupoid* G *with* $G^2$ *isomorphic to* U *is an inflation of* $G^2$.

In Theorem 2, if U is a semigroup then it is a semilattice of groups. Hence, Theorem 2 generalises Tamura's result that if S is a semigroup and $S^2$ is a semilattice of groups then S is an inflation of $S^2$ [18]. Theorem 2 also implies that a right modular generalised inflation of an idempotent, right modular groupoid U is an inflation of U. As noted above, a (semigroup) generalised inflation S of its subsemigroup U, where U is a union of groups, need ***not*** be an inflation of U.

Note that any semigroup S in which $S^2$ has an identity element is an inflation of $S^2$. Consequently, any generalised inflation S of a subsemigroup U of the semigroup S, in which $U^2$ has an identity element, is an inflation of $U^2 = S^2$. This result holds because in the definition of a generalised inflation S of its subsemigroup U, S is a semigroup. However, a generalised inflation G of a semigroup U need not be a semigroup if we treat U as a groupoid, since the definition of a generalised inflation G of its subgroupoid U does not require G to be a semigroup. For example, consider the two-element chain $C = \{a,b\}$, where $ab = ba = b$, $a = a^2$ and $b = b^2$. Let $U = C^1$, $G = U \cup \{x\}$ ($x \notin U$), $1\alpha_x = 1 = \beta_x 1$, $a\alpha_x = b = b\alpha_x$, $\beta_x a = a$ and $\beta_x b = b$.

Then G is a generalised inflation of U, with $x \in G_1$, but $x^2 a = 1a = a \neq b = b^2 = b \cdot ba = b\alpha_x \cdot [(a\alpha_x)a] = x \cdot xa$ and so G is not a semigroup. Note that U is a right modular groupoid with an identity element. However, the groupoid G is not right modular. The multiplication table for $G = U \cup \{x\}$ is as follows.

**Example 2.**

| G | a | b | 1 | x |
|---|---|---|---|---|
| a | a | b | a | b |
| b | b | b | b | b |
| 1 | a | b | 1 | 1 |
| x | b | b | 1 | 1 |

This leaves us with the following natural question**:** If G is a right modular generalised inflation of its subgroupoid U and U has an identity element, then is G an inflation of U *?*

**Theorem 3.** *If* G *is a right modular generalised inflation of its subgroupoid* U *and* U *has a left identity element, 1 say, then* G *is an inflation of* U.

*Proof:* Let $\{x, y\} \subseteq G$. Since G is a generalised inflation of U, $G^2 \subseteq U$. Then, using Result 1, $xy = (1 \cdot 1) \cdot xy = 1x \cdot 1y$. Now, $G = \bigcup G_u (u \in U)$, where for $u \in U$, $G_u = \{a \in G : 1a = u\}$, is a disjoint union satisfying $u \in G_u (u \in U)$. Also, for $x \in G_u$ and $y \in G_v$ ($\{u,v\} \subseteq U$) we have $xy = uv$. By definition then, G is an inflation of U. This completes the proof of Theorem 3.

**Theorem 4.** *If $G$ is a right modular groupoid and $x \notin G$ then $G \cup \{x\}$ is a right modular generalised inflation of $G$ if and only if there exist right and left mappings, $\alpha_x$ and $\beta_x$ respectively, on $G$ and an element $c \in G$ such that for all $\{a,b\} \subseteq G$, (1) $(a\alpha_x a)b = (ba)\beta_x(ba)$, (2) $(a\beta_x a)b = (b\beta_x b)a$ and (3) $(c\alpha_x \cdot \beta_x c) \cdot a = (a\beta_x a) \cdot \beta_x (a\beta_x a)$.*

*Proof:* ($\Rightarrow$) Since $G \cup \{x\}$ is a right modular generalised inflation of $G$, the mappings $\alpha_x$ and $\beta_x$ exist. Let $x \in G_c$. The equations $xa \cdot b = ba \cdot x$, $ax \cdot b = bx \cdot a$ and $x^2 \cdot a = ax \cdot x$ are valid, since $G \cup \{x\}$ is a right modular groupoid. By the definition of a generalised inflation of $G$, (1), (2) and (3) follow.

($\Leftarrow$) Define $G_a = a (a \neq c)$ and $G_c = \{c, x\}$. Then $G \cup \{x\} = \bigcup G_a (a \in G)$, with $G_a \cap G_b = \emptyset$ if $a \neq b$. Define $d\alpha_c = c = \beta_c d$ for any $\{c,d\} \subseteq G$. Then since (1), (2) and (3) are valid, we can use the $\alpha$'s and $\beta$'s to define a generalised inflation in which $xa \cdot b = ba \cdot x$, $ax \cdot b = bx \cdot a$ and $x^2 \cdot a = ax \cdot x$. Such a generalised inflation satisfies the equation $yz \cdot w = wz \cdot y$ and is therefore a right modular generalised inflation of $G$. This completes the proof of Theorem 4.

It is straightforward to show that the four inflations $U \cup \{x\}$ of $U$, where $U$ is as in Example 1, are isomorphic. So of the 262,144 groupoids on $G = U \cup \{x\}$ in which $G^2 \subseteq U$, and therefore in which $G$ may be a generalised inflation of $U$, there is (to within isomorphism) only one right modular generalised inflation of $U$. Clearly then, it is quite a powerful condition to require a generalised inflation to be right modular.

**Open Question.** The author does not know if a semigroup $G$ can be a generalised inflation of its commutative subsemigroup $U$ without being an inflation of $U$. Similarly, can a right modular groupoid $G$ be a generalised inflation of its (right modular) subgroupoid $U$ without being an inflation of $U$?

## References

[1] Amjid, F., Amjid, V., and Khan, M.: On some classes of Abel-Grassmann's groupoids, J. Adv. Res. Pure Math., **3,** (2011) 109-119.

[2] Clarke, G. T. and Monzo, R. A. R., A Generalisation of the Concept of an Inflation of a Semigroup, Semigroup Forum, **60** (2000) 172-186.

[3] Holgate, P.: Groupoids satisfying a simple invertive law. Math. Student **61**, (1992) 101-106.

[4] Jezek, J., Kepka, T.: Modular groupoids, Czech. Math. J., **34**, (1984) 477-487.

[5] Kamran, M. S., Mushtaq, Q., Finite AG-groupoid with left identity and left zero, IJMMS, **27(6)**, (2001) 387-389.

[6] Kazim, M.A., Naseeruddin, M.: On almost semigroups, .Portugaliae Mathematica, **36**, (1977) 41-47.

[7] Khan, M., Mushtaq, Q.: Ideals in left almost semigroups, arXiv:0904.1635v1 [math. GR] (2009).

[8] Khan, M., Mushtaq, Q.: Topological structure on Abel-Grassmann's groupoids, arXiv:0904.1650v1 [math.GR].

[9] Khan,S., Khan, M.: On decomposition of Abel-Grassmann's groupoids, Semigroup Forum, (submitted).

[10] Monzo, R. A. R.: Further results in the theory of generalised inflations of semigroups. Semigroup Forum, **76** (2008) 540-560.

[11] Monzo, R.A.R.: On the structure of Abel Grassmann Unions of groups, Int. J. of Algebra, **26**, (2010) 1261-1275.

[12] Monzo, R.A.R.: Power groupoids and inclusion classes of Abel Grassmann groupoids Int. J. of Algebra, **5(2)**, (2011) 91-105.

[13] Mushtaq, Q., Iqbal, M.: Partial ordering and congruences on LA-semigroups. Indian J. pure appl. Math., **22(4)**, (1991) 331-336.

[14] Mushtaq, Q.: Zeroids idempoids in AG-groupoids, Quasigroups and Related Systems, **11**, (2004) 79-84.

[15] Protic, P.V., Stevanovic, N.: Abel-Grassmann's bands. Quasigroups and Related Systems **11,** (2004) 95-101.

[16] Protic, P.V., Stevanovic, N.: Inflations of the AG groupoids, Novi Sad J. Math., **29(1)**, (1999) 19-26.

[17] Protic, P.V., Stevanovic, N.: Composition of Abel Grassman's 3-bands. Novi Sad J. Math., **34 (2)**, (2004) 5-182.

[18] Tamura, T.: Semigroups Satisfying Identity $xy = f(x, y)$, Pacific J. Math. **31** (1962) 513-521.

[19] Wang, Q. and Wismath, S. L.: Generalised Inflations and Null Extensions, Discussiones Mathematicae, General Algebra and Applications, **24** (2004) 225-249.